\documentclass[alpha-refs]{wiley-article}
\usepackage{color}
\usepackage{ulem}
\usepackage{comment}
\usepackage{algpseudocode}
\usepackage{algorithmicx}
\usepackage[section]{algorithm}
\renewcommand{\algorithmicrequire}{\textbf{Input:}}
\renewcommand{\algorithmicensure}{\textbf{Output:}}
\usepackage{subcaption}
\usepackage{bm}  % For bold lower-case greek letters

\newtheorem{prop}{Proposition}

\newcommand{\bA}{\mathbf{A}}
\newcommand{\bB}{\mathbf{B}}
\newcommand{\bC}{\mathbf{C}}
\newcommand{\bD}{\mathbf{D}}
\newcommand{\bd}{\mathbf{d}}
\newcommand{\bE}{\mathbf{E}}
\newcommand{\bF}{\mathbf{F}}
\newcommand{\bff}{\mathbf{f}}
\newcommand{\bG}{\mathbf{G}}
\newcommand{\bI}{\mathbf{I}}
\newcommand{\bK}{\mathbf{K}}

\newcommand{\bP}{\mathbf{P}}
\newcommand{\bQ}{\mathbf{Q}}
\newcommand{\bR}{\mathbf{R}}

\newcommand{\bS}{\mathbf{S}}
\newcommand{\bT}{\mathbf{T}}
\newcommand{\bU}{\mathbf{U}}
\newcommand{\bUb}{\mathbf{U}_\mathbf{B}}
\newcommand{\bUbj}{\mathbf{U}_{\mathbf{B}_j}}
\newcommand{\bUR}{\mathbf{U}_\mathbf{R}}
\newcommand{\bV}{\mathbf{V}}
\newcommand{\bW}{\mathbf{W}}
\newcommand{\bx}{\mathbf{x}}
\newcommand{\bxb}{\mathbf{x}^{\text{b}}}
\newcommand{\dx}{\delta \mathbf{x}}
\newcommand{\bY}{\mathbf{Y}}
\newcommand{\by}{\mathbf{y}}
\newcommand{\dz}{\delta \mathbf{z}}

\newcommand{\bDelta}{\mathbf{\Delta}}
\newcommand{\boeta}{\boldsymbol{\eta}}
\newcommand{\beps}{\boldsymbol{\epsilon}}
\newcommand{\bgamma}{\boldsymbol{\gamma}}
\newcommand{\bGamma}{\mathbf{\Gamma}}
\newcommand{\bLambda}{\mathbf{\Lambda}}
\newcommand{\bSigma}{\mathbf{\Sigma}}
\newcommand{\bomega}{\boldsymbol{\omega}}
\newcommand{\bOmega}{\mathbf{\Omega}}
\newcommand{\bmu}{\boldsymbol{\mu}}
\newcommand{\bpsi}{\boldsymbol{\psi}}
\newcommand{\bPsi}{\mathbf{\Psi}}
\newcommand{\bPhi}{\mathbf{\Phi}}
\newcommand{\bxi}{\boldsymbol{\xi}}

\usepackage{siunitx}

\papertype{Research Article}
\title{Accelerating an ensemble of variational data assimilations with randomized preconditioning}

\author[1,2]{Ieva Dau\v{z}ickait\.{e}}
\author[1,2]{Selime G\"{u}rol}
\author[3]{Mayeul Destouches}
\author[3]{Lo\"{i}k Berre}
\author[1,2]{Anthony T. Weaver}

\affil[1]{Cerfacs, Toulouse, France}
\affil[2]{CECI, Université de Toulouse, CERFACS/CNRS/IRD, Toulouse, France}
\affil[3]{Météo-France, CNRS, Univ. Toulouse, CNRM, Toulouse, France}

\corraddress{Ieva Dau\v{z}ickait\.{e}, CERFACS, 42 Av. Gaspard Coriolis, 31100 Toulouse, France. }
\corremail{dauzickaite@cerfacs.fr}

\fundinginfo{This project has received financial support from the French national programme LEFE (Les Enveloppes Fluides et l'Environnement).}

\runningauthor{Dau\v{z}ickait\.{e} et al.}

\begin{document}

\maketitle

\begin{abstract}
Ensembles of variational data assimilations (EDA) require solving systems of linear equations with iterative methods.
The solution process can be accelerated using a limited memory preconditioner constructed with approximations of the leading eigenpairs of the Hessian matrix. Randomized methods for low-rank matrix approximations provide a feasible approach for computing these eigenpairs. These methods use a random sketching matrix to obtain a low-rank representation of the Hessian matrix, which is then used for computing the eigendecomposition. The sketching matrix highly influences the quality of the approximation. In this paper, we show how the structure of the EDA can be exploited to construct a suitable sketching matrix, i.e., using the differences of the right-hand sides of the linear systems of equations. Idealised numerical experiments with the Lorenz-96 model show that the resulting preconditioner is able to accelerate the EDA solution process for all ensemble members, even if constructed from the control member only.
% Please include a maximum of seven keywords
\end{abstract}

\section{Introduction}

In numerical weather prediction, data assimilation (DA) combines observations of the meteorological state with a previous short-term forecast (background) while accounting for their respective uncertainties. The resulting state is used to initialize the next forecast. Ensemble approaches to data assimilation can provide multiple initial states for an ensemble forecast and can be also used to estimate flow-dependent error covariances (e.g., \citealp{hamill2000hybrid}, \citealp{lorenc2003potential}, \citealp{buehner2005ensemblederived}).
We focus on a specific ensemble approach, namely, an ensemble of data assimilations (EDA), which is used at some operational centres, including M\'{e}t\'{e}o-France and the European Centre for Medium-Range Weather Forecasts (ECMWF). The EDA consists of solving perturbed variational DA problems, 
which are constructed using perturbed inputs (observations and background); see, for example, \cite{pereira2006use}, \cite{bonavita2012use} and \cite{bouyssel2020}. 
Running an EDA is made computationally feasible by using approximations, for example, considering a limited number of ensemble members and integrating the model at lower resolution than that of the deterministic forecast model. 

The nonlinear variational DA problems are generally solved using the incremental approach \citep{courtier1994strategy}, which is a variant of the Gauss-Newton method \citep{gratton2007TGN, nocedal2006numerical}. Each iteration of this approach involves minimising a quadratic cost function obtained by linearising the nonlinear operators. In the EDA setting, it is common to use only one iteration of this incremental approach to control the computational cost \citep[e.g.,][]{bouyssel2020}. In this ensemble setting, this iteration requires minimising a quadratic cost function corresponding to a control unperturbed member and an ensemble of quadratic cost functions defined with the perturbed inputs. The minimisations are performed in parallel, via iterative solution of the associated systems of linear equations. 
The EDA structure can be exploited to accelerate the minimization. For instance, \cite{desroziers2012accelerating} use an iterative method for minimizing the quadratic cost functions and show how the recurrences of the solver can be used to construct an initial guess of the solution for related minimizations with closely related Hessian matrices. This approach can also be used to generate a preconditioner for all the ensemble problems, however it implies that part of the ensemble problems must be solved in sequence rather than in parallel, which may not always be possible in an operational setting.    
\cite{mercier2018block,mercier2019speeding} consider using block Krylov subspace methods to solve all the linear systems simultaneously. This approach considerably reduces the required wall-clock time when operating in observation space, that is, solving the dual problem, due to the reduced cost of orthogonalization and communication compared to the primal approach. However, block Krylov subspace methods rely on the assumption that all the linear systems are formulated using the same Hessian matrix, which only holds when the same linearisation state is used for every ensemble member.

In this paper, we focus on preconditioning as a mean to accelerate the minimisations of the EDA. 
Limited memory preconditioners (LMPs) are used operationally for the incremental strong-constraint four-dimensional variational problem (4D-Var) \citep{tshimanga2008limited}. The LMP variant constructed with approximate eigenpairs, called the spectral LMP, was originally proposed by \cite{fisher1998minimization}. 
In an operational context, each iteration of the incremental approach (called outer iterations hereafter to avoid confusion with the inner iterations of an iterative Krylov subspace method) can provide at negligible extra cost such approximate eigenpairs of the Hessian matrix; this is done using recurrences of the iterative solver used to solve the linear system.
These approximate eigenpairs can then be used for the subsequent outer iterations, starting from the second one. 
However, this is not applicable in an EDA setting with a single outer iteration. 
In this case, a feasible strategy for constructing a spectral LMP is to employ randomized methods to build a low-rank matrix approximation (for example, \cite{halko2011finding} and \cite{martinsson2020randomized}). These methods can be considered as random dimension reduction techniques, where a large matrix is approximated by a significantly smaller matrix, called the sketch. The sketch is obtained by computing a matrix-matrix product of the large matrix and a random sketching matrix. The sketch can be used to obtain an approximate solution (sketch-and-solve approach) or to construct a preconditioner for an iterative method (sketch-and-precondition approach). 

Randomized methods have been explored to accelerate the solution of systems of linear equations in variational DA when multiple outer iterations of the incremental approach are performed to solve a single DA problem. \cite{bousserez2018optimal} and \cite{bousserez2020enhanced} propose and analyse iterative methods that are based on the sketch-and-solve approach and are designed to replace the deterministic iterative solver for 4D-Var.
These methods are highly parallelizable and can reduce the minimization wall-clock time at a price of increased energy consumption when enough computational cores are available. However, the quality of the obtained solution is highly sensitive to the distribution of the eigenvalues of the Hessian matrix of the quadratic cost function. The sketch-and-precondition approach for variational DA using spectral LMPs has been studied with toy dynamical models by \cite{dauvzickaite2021randomised}, \cite{scotto_phdthesis_2022}, and \cite{subrahmanya2024randomized}. \cite{dauvzickaite2021randomised} show that spectral LMPs constructed using randomized eigendecomposition methods can accelerate the solution of an incremental weak-constraint 4D-Var problem. \cite{scotto_phdthesis_2022} proposes randomized algorithms for computing a generalized eigendecomposition to construct spectral LMPs when solving the incremental problem, using a variant of the preconditioned conjugate gradient method proposed by~\citep{DerberRosati89} that incorporates the spectral LMP preconditioner \citep{gurol2013solving}. \cite{scotto_phdthesis_2022} also investigates the use of randomized algorithms for the dual problem, and illustrates their effectiveness for the incremental strong-constraint 4D-Var problem. \cite{subrahmanya2024randomized} investigate an adaptive approach to constructing randomized spectral LMPs for the incremental strong-constraint 4D-Var problem, although the number of required sequential matrix-vector products may be not feasible in an operational setting.

In this paper, we investigate how the sketch-and-precondition approach can be used for an EDA. The sketching step requires computing a matrix-matrix product with the Hessian matrix, which involves running the tangent linear and adjoint models a number of times. However, this can be done in parallel and thus the required wall-clock time can be equivalent to one run of the tangent linear and adjoint models. While the quality of the approximation can be improved by using the power method, that is, computing a few sequential matrix-matrix products with the Hessian matrix, the cost of this approach may be prohibitive in an operational setting. We show that the structure of the EDA, namely, the right-hand side vectors of the perturbed linear systems of equations, can be exploited in order to generate a suitable sketching matrix. 
Under simplifying assumptions, the proposed approach is equivalent to sampling the random sketching matrix from a Gaussian distribution whose covariance matrix incorporates some Hessian matrix information; such approaches have been investigated for randomized approximations of integral operators by \cite{boullegeneralization}, \cite{persson2025randomized}, and for matrices arising in DA by \cite{scotto_phdthesis_2022} as we further describe in Section~\ref{sec:sketching_matrices}.
When the randomized method uses the proposed sketching matrix, we can obtain a higher quality approximate eigendecomposition, which in turn gives a more effective spectral LMP. Such a spectral LMP is shown to be efficient in preconditioning both the control and the perturbed ensemble problems. Our main contributions thus are:
\begin{itemize}
    \item investigating the use of spectral LMP generated with randomized methods in an EDA setting;
    \item proposing the use of right-hand side vectors of the linear systems of equations for constructing a sketching matrix and analysing its properties.
\end{itemize}

The structure of this paper is as follows. We start by briefly introducing the EDA framework in section~\ref{sec:EDA}. The spectral LMPs, randomized methods, and the new sketching matrix are discussed in section~\ref{sec:precond_random_methods}. In section~\ref{sec:numerical_experiments}, numerical examples with the Lorenz-96 model illustrate the sensitivities and the potential of the proposed preconditioning approach. We conclude in section~\ref{sec:conclusions}.

\section{Ensemble Data Assimilation: problem formulation}\label{sec:EDA}

\subsection{Incremental 4D-Var}
In strong-constraint 4D-Var, we consider an assimilation window from time $t^{(0)}$ to $t^{(N)}$, and search for the model state at the initial time $\bx \in \mathbb{R}^n$. The optimal state $\bx_* \in \mathbb{R}^n$ is defined as the minimizer of the following cost function:
\begin{equation}\label{eq:4dvar_nonlinear}
    \mathcal{J}(\bx) = \frac{1}{2} ( \bx - \bxb )^{\rm T} \bB^{-1} ( \bx - \bxb ) + \frac{1}{2} \left( \by - \mathcal{G}(\bx) \right)^{\rm T} \bR^{-1} \left( \by - \mathcal{G}(\bx) \right),
\end{equation}
where $\bxb \in \mathbb{R}^n$ is the background state, observations are collected into a vector $\by \in \mathbb{R}^p$, $\bB \in \mathbb{R}^{n \times n}$ and $\bR \in \mathbb{R}^{p \times p}$ are, respectively, the background- and observation-error covariance matrices, and the estimated state $\bx$ is mapped to the observation space at the correct time via the generalized nonlinear observation operator $\mathcal{G}$.
The computational time constraints in operational settings allow only a few integrations of the full non-linear model needed to evaluate $\mathcal{J}(\bx)$, thus optimization methods are generally stopped before full convergence.

The incremental approach, introduced to atmospheric DA by \cite{courtier1994strategy}, approximates $\bx_*$ via the truncated Gauss-Newton (GN) method \citep{gratton2007TGN,nocedal2006numerical}. In this paper, we consider a single outer iteration of the incremental approach. 
Here, we minimize the following quadratic cost function
\begin{equation}\label{eq:quadratic_cost_function}
    J(\dx) = \frac{1}{2} \dx^{\rm T} \bB^{-1} \dx + \frac{1}{2} (\bG \dx - \bd )^{\rm T}\bR^{-1}(\bG \dx - \bd ),
\end{equation}
where $\bG$ denotes the Jacobian matrix of $\mathcal{G}$ linearized around the trajectory of $\bxb$ and $\bd = \by - \mathcal{G}(\bxb)$ is the innovation vector.
The state estimate is updated as $\bxb + \delta \tilde{\bx}_*$, where $\delta \tilde{\bx}_*$ is the approximate solution of \eqref{eq:quadratic_cost_function}. The linearized model $\bG$ can be run at lower resolution than the non-linear model $\mathcal{G}$, which makes the minimization tractable in operational settings.

\subsection{Ensemble formulation}

The EDA is obtained by considering \eqref{eq:4dvar_nonlinear} with independently perturbed inputs. In the following, we use a subscript $j$ with $j=1,\dots,L$ to denote the inputs for a $j$th perturbed EDA member. 
The observations are perturbed explicitly, that is for a $j$th ensemble member, we compute
\begin{equation}
    \by_j = \by + \bUR \boeta_j^o,
\end{equation}
where $\bR = \bUR \bUR^\text{T}$ and $\boeta_j^o \sim \mathcal{N}(\mathbf{0},\mathbf{I})$, so that $\by_j\sim\mathcal{N}(\by,\bR)$. 

The perturbed background state $\bxb_j$ comes from a previous ensemble forecast, which uses the results of a previous EDA cycle. This way the background states are affected by the previously perturbed observations and the subsequent ensemble forecast. The background states are also affected by the model error: if the EDA system does not explicitly account for it in ensemble forecasts, it can be simulated in $\bxb_j$ using, for example, inflation \citep{raynaud2012accounting}. Resulting background perturbations $\beps^\text{b}_j$ defined by departures from the control unperturbed background, namely $\beps^\text{b}_j = \bxb_j - \bxb$, can then be considered as random draws of $\bB$.

We consider the case where $\bG_j$ is linearized around $\bxb_j$. Other approaches where the non-linear operators are linearized around the same state could also be used.
For example, the non-linear operators are linearized around the ensemble background mean when using block methods to solve the perturbed problems in the EDA \citep{mercier2018block,mercier2019speeding}. 
Since the specified operators $\bG_j$ and $\bB_j$ may depend on the linearization state (e.g., \citealp{geer2011observation} for $\bG_j$ and \citealp{fisher2003background} for $\bB_j$), the error covariance matrices $\bR_j$ (as it accounts for error in $\bG_j$) and $\bB_j$ may be different for separate ensemble members when using different linearization states. 
Thus, the quadratic cost function for the $j$th ensemble member is given by \eqref{eq:quadratic_cost_function} with the quantities $\dx_j$, $\bxb_j$, $\bd_j$, $\bG_j$, $\bR_j$ and $\bB_j$ dependent on the $j$th member. The unperturbed problem \eqref{eq:quadratic_cost_function} is called the control member.

\subsection{Solving linear systems of equations in the EDA}\label{sec:solving_linear_systems}

The quadratic cost function \eqref{eq:quadratic_cost_function} is minimized by solving the following linear system of equations:
\begin{equation}\label{eq:control_system}
    \left( \bB^{-1} + \bG^\text{T} \bR^{-1} \bG \right) \dx = \bG^\text{T} \bR^{-1} \bd,
\end{equation}
where $\bB^{-1} + \bG^\text{T} \bR^{-1} \bG$ is the Hessian matrix. An equivalent system has to be solved for each perturbed ensemble member. 
Equation \eqref{eq:control_system} must be solved iteratively due to the large size of the problem and the fact that the matrices are generally only available through matrix-vector products, without explicit access to their entries.
Note that the Hessian is symmetric positive definite. The conjugate gradient (CG) method or its Lanczos variant (Lanczos-CG) are widely used for solving such systems \citep{hestenes1952methods,frommer1999fast,saad2003iterative}.
These solvers are also guaranteed to reduce the value of \eqref{eq:quadratic_cost_function} on every iteration. 
Performance of CG is sensitive to the distribution of the eigenvalues of the Hessian matrix, as well as the right-hand side vector though this latter influence is often disregarded \citep{liesen2013krylov}. In general, clusters of eigenvalues close to zero can cause stagnation of the method for a number of iterations, and it may be beneficial to have clusters of eigenvalues away from zero; see, for example, \cite{carson2024towards} for a deeper discussion on CG convergence.
Preconditioning techniques, e.g., \cite{wathen2015preconditioning}, can be used to map the given problem to another problem, which has a more favourable eigenvalue distribution and allows the solution to the original problem to be easily recovered.

The Hessian matrix in Eq.~\eqref{eq:control_system} includes the background error covariance matrix $\bB$, which is usually ill-conditioned, that is, the ratio between its largest and smallest eigenvalues is large \citep{lorenc1997development}. This severely affects the performance of the iterative solver without preconditioning. The modelling of $\bB$ as $\bUb\bUb^\text{T}$, where $\operatorname{rank}(\bUb) =m \leq n$, can be used to introduce the change of variable for first-level preconditioning, i.e., 
\begin{equation}\label{eq:first_level_prec}
    \dx = \bUb \dz.
\end{equation}
Then \eqref{eq:control_system} can be written in terms of the split-preconditioner as \citep{fisher1998minimization,menetrier2015overlooked,saad2003iterative}
\begin{equation}\label{eq:system_with_CVT}
    \left( \bI + \underbrace{\bUb^\text{T} \bG^\text{T} \bR^{-1} \bG \bUb}_{\bA} \right) \dz = \bUb^\text{T} \bG^\text{T} \bR^{-1} \bd,
\end{equation}
where $\dx$ is retrieved from $\dz$ using \eqref{eq:first_level_prec}.
The preconditioned Hessian in \eqref{eq:system_with_CVT} is expressed as a sum of an identity matrix and a low-rank symmetric positive semidefinite matrix $\bA$. 
Thus $\bI + \bA$ has a cluster of eigenvalues equal to one, and the rest are equal to $1+\lambda_i(\bA)$, where $\lambda_i(\bA)$, $i=1,\dots,\operatorname{rank}(\bA)$, are the non-zero eigenvalues of $\bA$. The eigenvectors of $\bA$ are also the eigenvectors of $\bI+\bA$. An additional preconditioner, a so-called second-level preconditioner, that tackles the largest eigenvalues may further improve the performance of iterative solvers.

\section{Preconditioning with randomized algorithms}\label{sec:precond_random_methods}
With the second-level preconditioning, we aim to approximate the inverse of $\bI+\bA$. 
Since in operational settings only a limited amount of computational resources can be used to perform matrix-vector products with a preconditioner, a widely adopted preconditioner is the so-called inexact spectral LMP built with Ritz pairs, an approximation of the eigenpairs of $\bI+\bA$ \citep{fisher2009data,gratton2011class}. These Ritz pairs are computed from the previous outer iteration of the Gauss-Newton method, which makes this approach not applicable for EDAs with a single outer iteration. Alternatively,
we consider using parallelizable randomized methods to estimate eigenpairs. A scaled version of the spectral LMP as well as randomized methods for computing a truncated eigendecomposition are discussed in the rest of this section.

\subsection{Spectral limited memory preconditioner}
 Let $\bA = \bS \bLambda \bS^\text{T}$ be an eigendecomposition of $\bA$, and consider a partition $\bLambda = \operatorname{diag}(\bLambda_k, \bar{\bLambda}_k)$, where $\bLambda_k = \operatorname{diag}(\lambda_1(\bA),\dots, \lambda_k(\bA))$  and $ \bar{\bLambda}_k = \operatorname{diag}(\lambda_{k+1}(\bA),\dots, \lambda_n(\bA))$ with eigenvalue ordering $\lambda_1(\bA) \geq \lambda_2(\bA) \geq \dots \geq \lambda_n(\bA) \geq 0$. Similarly, we partition $\bS$ as $\bS = \begin{bmatrix} \bS_k & \bar{\bS}_k \end{bmatrix}$, where $\bS_k$ is of size $n \times k$ and $\bar{\bS}_k$ is of size $n \times (n-k)$. Then the scaled spectral LMP for preconditioning $\bI+\bA$ is defined as
\begin{equation}
    \bP_{\theta} = \bI + \bS_k \left( \theta (\bLambda_k + \bI)^{-1} - \bI \right) \bS_k^\text{T},
\end{equation}
where $\theta >0$ is a scaling factor. $\bP$ can be factorized as $\bP_{\theta} = \bU_{\theta} \bU_{\theta}^\text{T}$ with
\begin{equation}\label{eq:lmp}
    \bU_{\theta} = \bU_{\theta}^\text{T} = \bI + \bS_k \left( \sqrt{\theta} (\bLambda_k +\bI)^{-1/2} - \bI \right) \bS_k^\text{T}.
\end{equation}
The preconditioner $\bP_{\theta}$ only affects the eigenvalues which are used in its construction. That is, since $\bP_{\theta}$ is constructed with the $k$ leading eigenpairs of $\bI + \bA$, the preconditioned matrix $\bU_{\theta}(\bI + \bA)\bU_{\theta}^\text{T}$ have these leading eigenvalues moved to a cluster at $\theta$ and the other eigenvalues remain unchanged; see \cite{gratton2011class} for a general analysis of LMPs. An appropriately chosen $\theta$ can thus give both a smaller condition number and a clustering of the eigenvalues better suited to the CG/Lanczos-CG method. The scaling factor $\theta$ can also be chosen to influence the behaviour of the iterative method, especially in the first few iterations, which are of most importance in DA applications due to the early termination, as shown by \cite{diouane2024efficient}. They propose choosing $\theta = (\lambda_k(\bI+\bA) + \lambda_n(\bI+\bA))/2$, which ensures that the iterates of preconditioned CG (PCG) with the scaled spectral LMP closely match those of deflated CG \citep{saad2000deflated}, and discuss other choices for $\theta$. For further details, see \cite{diouane2024efficient}. 
Note that in our case, $\lambda_n(\bI+\bA)=1$, because $\bA$ is symmetric positive semidefinite and has no more than $p$ non-zero eigenvalues with $p<n$.

In practice, \eqref{eq:lmp} is constructed with Ritz pairs that approximate the leading eigenpairs of $\bI+\bA$. The Ritz pairs are obtained in the previous outer iteration and are used to construct \eqref{eq:lmp} only if they are of sufficient accuracy; an accuracy estimate for each Ritz pair being available cheaply as a by-product of the Ritz pair estimation. 
In this paper, we consider constructing the spectral LMP employing the randomized eigendecomposition of $\bI+\bA$.

\subsection{Randomized eigendecomposition}
Randomized algorithms for computing a low-rank approximation of a given matrix constitute a large area of research (see, for example, the surveys by \citealp{halko2011finding} or \citealp{martinsson2020randomized}). 
These algorithms employ randomness to reduce the dimension of the problem and use matrix-matrix multiplications, which can be implemented efficiently on modern high-performance computing systems. 

Algorithms for computing a randomized eigendecomposition (REVD) of an $n \times n$ matrix $\bA$ usually proceed in two stages. In the first stage, the range of $\bA$ is approximated by an $n \times \ell$ sketch $\bA \bOmega$ using a random sketching matrix $\bOmega$. The second stage is deterministic and consists of using the approximation of the range to map $\bA$ to a smaller matrix and computing its eigendecomposition (EVD). The quality of the approximation depends on the choice of the sketching matrix $\bOmega$, and the chosen projection from $\bA$ to the smaller matrix. Note that there is no computationally inexpensive way to evaluate the quality of the approximate eigenpairs coming from the randomized methods. We discuss specific choices of the sketching matrices in the next two subsections, and present an REVD algorithm in subsection~\ref{sec:nystrom}.

\subsubsection{Sketching matrices}\label{sec:sketching_matrices} 

Various sketching matrices can be used depending on the structure of the problem. In our case, we have to pass the sketching matrix $\bOmega = \begin{pmatrix}
    \bomega_1 & \bomega_2 & \dots & \bomega_{\ell}
\end{pmatrix}$ column-by-column to a procedure that returns $\bA \bomega_i$. Note that computations with each $\bomega_i$ can be performed independently in parallel. A well-studied sketching matrix is a matrix with random standard Gaussian entries. 

The quality of the approximation may be improved if $\bOmega$ is enriched with information on the subspace that we are aiming to approximate. The spectral LMP in \eqref{eq:lmp} requires approximations to the $k$ largest eigenvalues and the corresponding eigenvectors. An approximation to the eigenpair associated with the largest eigenvalue can be computed using the power method as $\bA^{q-1} \bff$ for some vector $\bff$ and $q>1$. Its generalisation, called the subspace iteration, can be employed for obtaining approximations of $k$ largest eigenpairs. In this case, one considers $\bA^{q-1} \bF$, where $\bF$ is of size $n \times \ell$ with $\ell \geq k$. We can thus employ this idea for generating the sketching matrix as
\begin{equation}
    \bOmega = \bA^{q-1} \bPsi,
\end{equation}
where $\bPsi = \begin{pmatrix}
    \bpsi_1 & \bpsi_2 & \dots & \bpsi_{\ell}
\end{pmatrix}$, the column vectors $\bpsi_i$ are independently sampled from standard Gaussian distribution and $q\geq1$ is an integer power parameter.

Another interpretation of the generalised power method can be obtained using the well-known fact that a linear transformation of a Gaussian random vector is a Gaussian random vector \citep[see, for example, Theorem 5.2 of][]{hardle2007applied}. Namely, if we have a random vector $\boldsymbol{X} \sim \mathcal{N}(\bmu,\boldsymbol{C}_X)$, where $\bmu$ is a vector of length $n$ and $\boldsymbol{C}_X$ is of size $n \times n$, then for an $m \times n$ matrix $\boldsymbol{T}$, the following holds $\boldsymbol{TX} \sim \mathcal{N}(\boldsymbol{T}\bmu, \boldsymbol{T} \boldsymbol{C}_X \boldsymbol{T}^\text{T})$. Thus, we can interpret using a sketching matrix $\bOmega = \bA^{q-1} \bPsi$ as sampling the columns of $\bOmega$ from a Gaussian distribution with zero mean and covariance matrix $\bA^{2(q-1)}$.

Note that computing $\bA^{q-1} \bpsi_i$ requires integrating the tangent linear model and adjoint model $q-1$ times sequentially. Setting $q$ to values larger than one or two may be prohibitive in an operational setting. It is desirable then to seek an approximation of the power method. \cite{di2024general} propose exploiting the structure of $\bA$ by using 
\begin{equation}\label{eq:omega_UbT}
    \bOmega = \bUb^\text{T} \bUb \bPsi \quad \operatorname{or} \quad \bOmega = \bUb^\text{T} \bPsi,
\end{equation}
where $\bPsi$ has an appropriate number of rows. Here we take into account that $\bUb$ may be rectangular, that is, 
$\dx$ and $\dz$ may be of different sizes, and thus the sizes of $\bOmega$ and $\bPsi$ have to be chosen accordingly.
Using the choices in \eqref{eq:omega_UbT} is equivalent to drawing $\bomega_i$ from a Gaussian distribution with covariance matrices $\bUb^\text{T}\bB\bUb$ and $\bB$, respectively; note that if $\bUb^\text{T}=\bUb$ then $\bUb^\text{T}\bB\bUb = \bB^2$. The proposed $\bOmega$ variants are in the range of $\bUb^\text{T}$ as is $\bA$ and can be obtained without the expensive integration of the tangent linear and adjoint models. 

In the next subsection, we further examine the structure of the EDA framework (that is, when solving perturbed linear systems), in order to obtain improved approximations of the generalized power method that incorporate information on both observations and model trajectories.

\subsubsection{Exploring the structure of the EDA}\label{sec:sketching_using_eda_structure}

Observe that randomness is already used in generating the ensemble. For each ensemble member, the right-hand side (RHS) vector in \eqref{eq:system_with_CVT} has to be computed before passing them to the iterative solvers. We can cheaply compute the difference between the RHS vectors of a $j$th ensemble member and the control member:
\begin{equation}\label{eq:gamma_def}
    \bgamma_j = \bUbj^\text{T} \bG^\text{T}_j \bR_j^{-1} \bd_j - \bU_{\bB}^\text{T} \bG^\text{T} \bR^{-1} \bd.
\end{equation}
Due to the randomly perturbed background and observations (and hence innovation vector $\bd_j$), $\bgamma_i$ is a random vector. Thus the matrix of RHS differences
\begin{equation}\label{eq:Gamma}
    \bGamma = \begin{pmatrix}
        \bgamma_1 & \bgamma_2 & \dots & \bgamma_{\ell}
    \end{pmatrix}
\end{equation}
can be used in the sketching step, that is, we can set $\bOmega = \bGamma$. In the following analysis, we assume that $\bGamma$ is full rank, which is an acceptable assumption in practice.
Note that $\ell$, i.e., the number of columns of $\bGamma$, is limited by the number of perturbed ensemble members $L$.
We analyse the properties of $\bGamma$ under simplifying assumptions in the following proposition.
\begin{prop}\label{prop:rhs_sketching_analysis}
Assume that $\bUbj^\text{T} \approx \bU_{\bB}^\text{T}$, $\bG^\text{T}_j \approx \bG^\text{T}$, and $\bR_j^{-1} \approx \bR^{-1}$, and 
that the same observations are used for the control and all ensemble members \footnote{This is the case for the current operational global EDA at Météo-France for instance. This assumption may not hold if, for example, the screening or variational quality control reject different observations for the ensemble members}. Denoting $\beps^\text{b}_j = \bxb_j - \bxb$, assume further that $\boeta^o_j$ and $\beps^\text{b}_j$ are independent, $\mathbb{E}\!\left[ \beps^\text{b}_j \right] = \mathbf{0}$, where $\mathbb{E}$ denotes the expectation operator, and $\operatorname{cov}(\beps^\text{b}_j,\beps^\text{b}_j) \approx \bB$. 
Then using a first order approximation $\mathcal{G}(\bxb_j) \approx \mathcal{G}(\bxb) + \bG \beps^\text{b}_j$, we can write
\begin{equation}\label{eq:gamma_approx}
    \bgamma_j \approx \bU_{\bB}^\text{T} \bG^\text{T} \bR^{-1} \left( \bU_{\bR} \boeta^o_j - \bG \beps^\text{b}_j \right), 
\end{equation}
which gives $\mathbb{E} \!\left[ \bgamma_j \right] = \mathbf{0}$ and
\begin{equation}\label{eq:cov_gamma}
    \operatorname{cov}(\bgamma_j,\bgamma_j) \approx  \bA + \bA^2.
\end{equation}
\end{prop}

\begin{proof}

The result in \eqref{eq:gamma_approx} follows from using $\bUbj^\text{T} \approx \bU_{\bB}^\text{T}$, $\bG^\text{T}_j \approx \bG^\text{T}$, and $\bR_j^{-1} \approx \bR^{-1}$ in \eqref{eq:gamma_def} and employing $\mathcal{G}(\bxb_j) \approx \mathcal{G}(\bxb) + \bG \beps^\text{b}_j$ to approximate $\bd_j = \by_j - \mathcal{G}(\bxb_j)$. 
Equation \eqref{eq:cov_gamma} is then obtained using the assumptions that $\boeta^o_j$ and $\beps^\text{b}_j$ are independent and have zero means, that is,
\begin{equation}
   \operatorname{cov}(\bgamma_j,\bgamma_j) \approx \mathbb{E}\! \left[ \bU_{\bB}^\text{T} \bG^\text{T} \bR^{-1} \left( \bU_{\bR} \boeta^o_j -  \bG \beps^\text{b}_j \right) \left( \bU_{\bB}^\text{T} \bG^\text{T}  \bR^{-1} \left( \bU_{\bR} \boeta^o_j -  \bG \beps^\text{b}_j \right)\right)^\text{T} \right] = \bA + \bA^2.
\end{equation}

\end{proof}

The results can be interpreted by considering the eigendecomposition $\bA = \bS \bLambda \bS^T$. Then $\bA + \bA^2 = \bS ( \bLambda + \bLambda^2) \bS^{\rm T}$. Thus, proposition~\ref{prop:rhs_sketching_analysis} establishes that the covariance of $\bGamma$ has the same eigenvectors as $\bA$ and gives more relative weight to the eigenvectors associated with the larger eigenvalues. As a consequence, we can expect that using $\bGamma$ as a sketching matrix may have a positive effect on the quality of the approximation.  

We can relax the assumption $\bG^\text{T}_j \approx \bG^\text{T}$ (i.e., that the different linearization states do not affect $\bG_j$) in Proposition~\ref{prop:rhs_sketching_analysis} and instead consider the case where $\bG_j = \bG +\bE_j$. Then using the rest of the assumptions in Proposition~\ref{prop:rhs_sketching_analysis}
and assuming that $\boeta^o_j$ is independent of the pair $\beps^\text{b}_j$ 
and $\bE_j$, we have
\begin{equation}\label{eq:cov_gamma_Ej_term}
    \operatorname{cov}(\bgamma_j,\bgamma_j) \approx \bA + \bA^2 + \bDelta,
\end{equation}
where $\bDelta$ is a perturbation in the range of $\bU_{\bB}^\text{T}$ and its size is dominated by $\Vert\mathbb{E}\! [\bE_j] \Vert_2$; see Appendix~\ref{app:perturbation_gamma} for the derivation.

\subsubsection{Nystr\"{o}m approximation}\label{sec:nystrom}
Since $\bA$ is symmetric positive semidefinite, the preferred approximation of $\bA$ is expressed in the following form called the Nystr\"{o}m approximation, see, e.g. Lemma 5.2 of \cite{tropp2023randomized},
\begin{equation}\label{eq:nystrom_approx}
    \bA_N \coloneqq \bA \bOmega \left( \bOmega^\text{T} \bA \bOmega \right)^{\dagger} (\bA \bOmega)^\text{T},
\end{equation}
where $\dagger$ denotes the Moore-Penrose pseudoinverse.
The approximation can be written as $\bA_N = \bK \bP_{\bK \bOmega} \bK$, where $\bK$ is the unique square-root of $\bA$ such that $\bK\bK = \bK^\text{T}\bK = \bA$, and $\bP_{\bK\bOmega}$ is an orthogonal projector onto the range $\bK \bOmega$. As discussed in the previous sections, the quality of the approximation depends on the range of $\bOmega$ with $\operatorname{range}(\bOmega) = \operatorname{range}(\bS_k)$ giving the smallest possible error. Deterministic and probabilistic bounds with non-standard Gaussian $\bOmega$ for different norms of $\bA - \bA_N$ are available in the literature, e.g., \cite{gittens2016revisiting}, \cite{persson2025randomized}. These bounds are in general pessimistic and we are unable to use them for guidance on the performance of the discussed choices of $\bOmega$.

Computing the approximation directly as in \eqref{eq:nystrom_approx} may cause numerical instabilities. We provide Algorithm~\ref{alg:nystrom} for computing $\bA_N$, and discuss other algorithmic variants in Appendix~\ref{app:Nystrom_algorithms}. Algorithm~\ref{alg:nystrom} uses oversampling, where the sketching matrix $\bPhi$ has more columns than the required rank of $\bA_N$. In this way, the sketch approximates $\bA$ in a larger subspace and the approximation quality can be improved. It is common to set the oversampling parameter $\ell - k$ to small values, e.g., 5 or 10.

In the case that computing matrix-vector products with $\bA$ is expensive, the leading cost of the randomized algorithm is $q$ matrix-matrix products. In our application, this is equivalent to $q$ sequential integrations of batches of $\ell$ independent instantiations of the tangent linear and adjoint models; each of the $\ell$ integrations can be done in parallel and thus the wall-clock time is equivalent to $q$ sequential integrations of the tangent linear and adjoint models.

\begin{algorithm}
\caption{$q$-pass Nystr\"{o}m approximation for symmetric positive semidefinite $\bA$}\label{alg:nystrom}
\algorithmicrequire \ symmetric positive semidefinite matrix $\bA \in \mathbb{R}^{n \times n}$, sketching matrix $\bPhi \in \mathbb{R}^{n \times \ell}$, a target rank $k \leq \ell$, number $q\geq1$ of matrix-matrix products with $\bA$ \\
\algorithmicensure \ diagonal matrix $ \bD_k \in \mathbb{R}^{k \times k}$ with approximations to the largest $k$ eigenvalues of $\bA$ on the diagonal, and $\hat{\bS}_k \in \mathbb{R}^{n \times k}$ with orthonormal columns approximating associated eigenvectors of $\bA$
\begin{algorithmic}[1]
\State Sketch $\bY = \bA \bPhi$ \label{algstep:sketch}
\State Compute the economic QR decomposition $\bY = \bQ_Y \bR_Y$
\For{$j = 2, \dots, q$} \Comment{power method}
        \State $\bPhi = \bQ_Y$
        \State $\bY = \bA \bPhi$
        \State Compute the economic QR decomposition $\bY = \bQ_Y \bR_Y$
\EndFor
\State  $\bW = \bPhi^\text{T} \bY$
\State Compute the upper triangular Cholesky factor: $\bW = \bC^\text{T} \bC$  
\State Compute via triangular solves $\bT = \bR_Y \bC^{-1}$ 
\State Compute the economic SVD $\bT =\bU_T \bSigma_T \bV_T^\text{T}$ \Comment{$\bV_T$ is not needed} 
\State Compute $\bS = \bQ_Y \bU_T$ and $\bD = \bSigma_T^2$
\State Remove oversampled columns/rows: $\hat{\bS}_k = \bS(:,1:k)$, $\bD_k = \bD(1:k,1:k)$
\end{algorithmic}
\end{algorithm}

\section{Numerical experiments}\label{sec:numerical_experiments}
We perform identical twin experiments in order to investigate the properties of the proposed randomized approximation approaches in an ensemble setting. 
An in-house implementation of a 4D-Var data assimilation system with the Lorenz-96 model is used \citep{lorenz1996predictability}. In this model, the dynamics for each $X^{(j)}$ in $\bx = \begin{pmatrix} X^{(1)} & X^{(2)} & \dots & X^{(n)}\end{pmatrix}^{\rm T}$ are determined by the following ODE
\begin{equation}
    \frac{dX^{(j)}}{dt} = - X^{(j-2)} X^{(j-1)} + X^{(j-1)} X^{(j+1)} - X^{(j)} + F,
\end{equation}
with $X^{(-1)} = X^{(n-1)}$, $X^{(0)} = X^{(n)}$, and $X^{(n+1)} = X^{(1)}$. The model simulates a motion of a single meteorological variable over $n$ points along a latitude circle. Setting $F=8$ ensures that the system shows chaotic behaviour. A fourth order Runge-Kutta scheme with a time step set to $2.5\times 10^{-2}$ is used. The true state $\bx^{\text{t}}$ is constructed by integrating the model for 5000 time steps starting with $\bx = \begin{pmatrix} 3 & 3 & \dots & 3\end{pmatrix}^{\rm T} + \bxi_1$, where $\bxi_1 \sim \mathcal{N}(\mathbf{0},\mathbf{I})$. We wish to investigate the performance of randomized methods with $k \ll n$, thus we choose $n=1.5 \times 10^3$.
The DA time window is equal to ten time steps, which corresponds to around 30 hours given the chosen time step size \citep{lorenz1996predictability}. This choice is in line with, for example, \cite{fertig2007comparative} who use a DA time window of 12h to 108h and \cite{zhang2009coupling} who consider a 24h to 60h time window, both in the context of 4D-Var experiments with the Lorenz-96 model.

The observation operator takes direct observations of fifty variables chosen uniformly every 30 grid points at three time steps ($t_0,t_3,t_7$) giving $p=1.5 \times 10^2$. The observations are perturbed using $\bR = \sigma_o^2 \bI$ with $\sigma_o = 5\times 10^{-2}$. The background is obtained as $\bxb =\bx^{\text{t}} + \bUb \bxi_2$ with the background-error covariance matrix $\bB=\bUb\bUb^{\rm T}$ set to $\bB = \sigma_b^2 \bC_\bB$, where $\sigma_b = 0.8$ and $\bC_\bB$ is a diffusion-based correlation operator \citep{goux2024impact}, and $\bxi_2 \sim \mathcal{N}(\mathbf{0},\mathbf{I})$. Correlation operators derived from diffusion operators are controlled by a smoothness parameter $M$, which corresponds to the number of diffusion iterations, and a length-scale parameter $D$, which is related to the square-root of the diffusion coefficient. A one-dimensional diffusion operator, as required in this study, represents an operator with an auto-regressive (AR) correlation function as its kernel. For large values of $M$ ($M \gtrsim 10$), the AR function is approximately Gaussian, while smaller values of $M$ correspond to correlation functions with fatter tails. The value of $D$ controls the spatial range of the correlations, but does not affect the roll-off characteristic of the correlation function. A comprehensive discussion is given in Section~2.3 of \cite{goux2024impact}. 
We test different values of the parameters $M$ and $D$ in the numerical experiments below (Section~\ref{sec:eig_sensitiv_B}).
The same $\bR$ and $\bB$ are used in the assimilation process and to generate the perturbed background and observations for the ensemble members, that is, $\bR_j=\bR$ and $\bB_j=\bB$. The perturbed backgrounds are set to $\bxb_j = \bxb+\bUb \boeta_j^\text{b}$, where $\boeta_j^\text{b} \sim \mathcal{N}(\mathbf{0},\mathbf{I})$. 
In each experiment, we generate $L = \ell = 20$ ensemble members. The generalized observations operators are linearized around the perturbed backgrounds.

We test different sketching matrices $\bOmega$ in \eqref{eq:nystrom_approx}, which correspond to different choices of the parameter $q$ and sketching matrix $\bPhi$ in Algorithm~\ref{alg:nystrom}. Namely,
\begin{itemize}
    \item $\bOmega = \bPsi$, where $\bPsi$ has independent standard Gaussian entries ($q=1$, $\bPhi = \bPsi$);
    \item $\bOmega = \bA \bPsi$ ($q=2$, $\bPhi = \bPsi$);
    \item $\bOmega = \bUb^\text{T} \bUb \bPsi= \bB \bPsi$ ($q=1$, $\bPhi = \bB \bPsi$);
    \item $\bOmega = \bUb^\text{T}\bPsi$  ($q=1$, $\bPhi = \bUb^\text{T} \bPsi$);
    \item $\bOmega = \bGamma$ as in \eqref{eq:Gamma} ($q=1$, $\bPhi = \bGamma$).
\end{itemize}
Note that $\bOmega = \bA \bPsi$ is the computationally most expensive variant, and we do not test higher values of $q$ as these are expected to be infeasible in operational settings. The variants $\bOmega = \bB \bPsi$ and $\bOmega = \bUb^\text{T}\bPsi$ have lower computational cost than $\bOmega = \bA \bPsi$ but they require $\ell$ additional products with $\bB$ and $\bUb^\text{T}$, respectively, compared to $\bOmega = \bPsi$; these matrix-vector products can be computed in parallel. The extra computations for obtaining $\bOmega = \bGamma$, which exploits the RHS vectors of the linear systems of equations in the ensemble setting, are negligible, but this approach requires communicating the RHS vectors of all the ensemble members before the iterative solvers are called.

We run Algorithm~\ref{alg:nystrom} to compute an approximation of $\bA$, that is, the Hessian of the control member. 
The output is then used to construct the scaled LMP in \eqref{eq:lmp}, which is used to precondition PCG when solving both control and perturbed ensemble members. The initial guess for PCG is set to a zero vector and the method is stopped after 20 iterations.

The code is written in \texttt{Python 3} and largely exploits the \texttt{SciPy} package \citep{2020SciPy-NMeth}.

\subsection{Results: eigenvalue sensitivity to $\bB$ }\label{sec:eig_sensitiv_B}
\begin{figure}[ht]
    \centering
    \captionsetup[subfigure]{justification=centering}
    \begin{subfigure}[t]{0.5\textwidth}
        \centering
        \includegraphics[width=\linewidth]{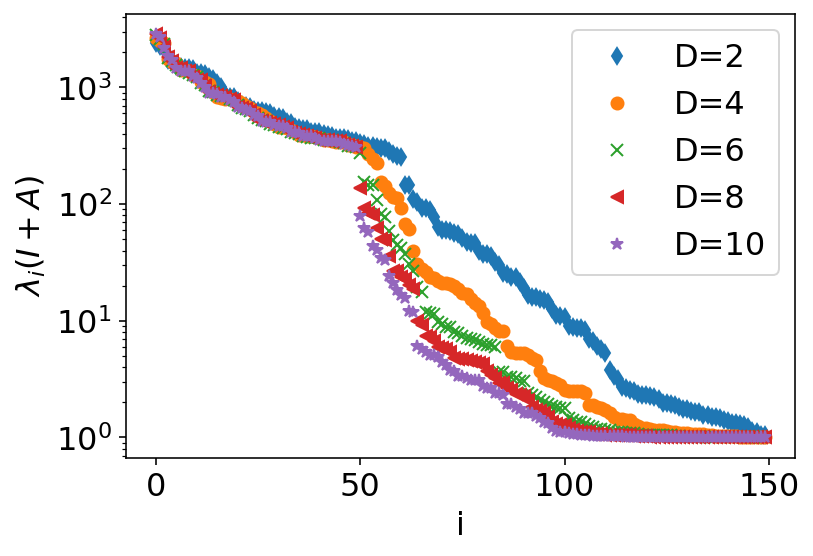}
        \caption{$M=10$; 150 eigenvalues}
    \end{subfigure}% 
    \begin{subfigure}[t]{0.5\textwidth}
        \centering
        \includegraphics[width=\linewidth]{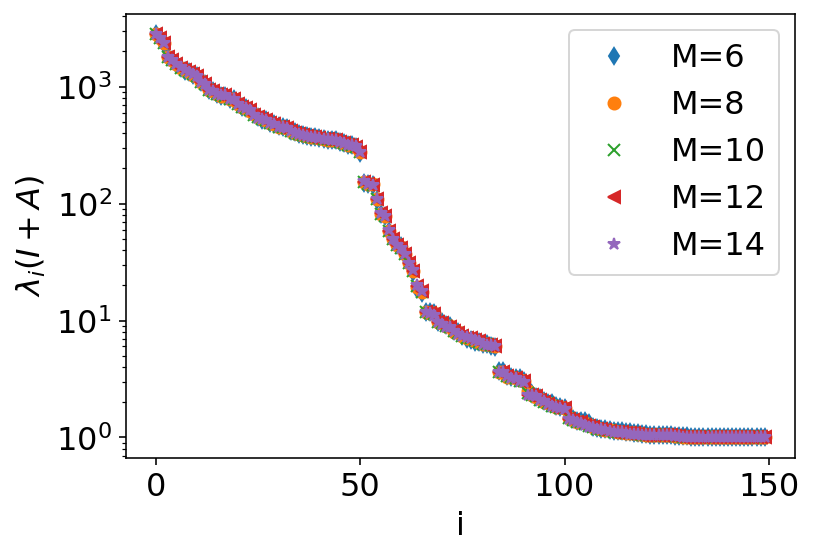}
        \caption{$D=6$; 150 eigenvalues}
    \end{subfigure}\vspace{0.5cm}
        \begin{subfigure}[t]{0.5\textwidth}
        \centering
        \includegraphics[width=\linewidth]{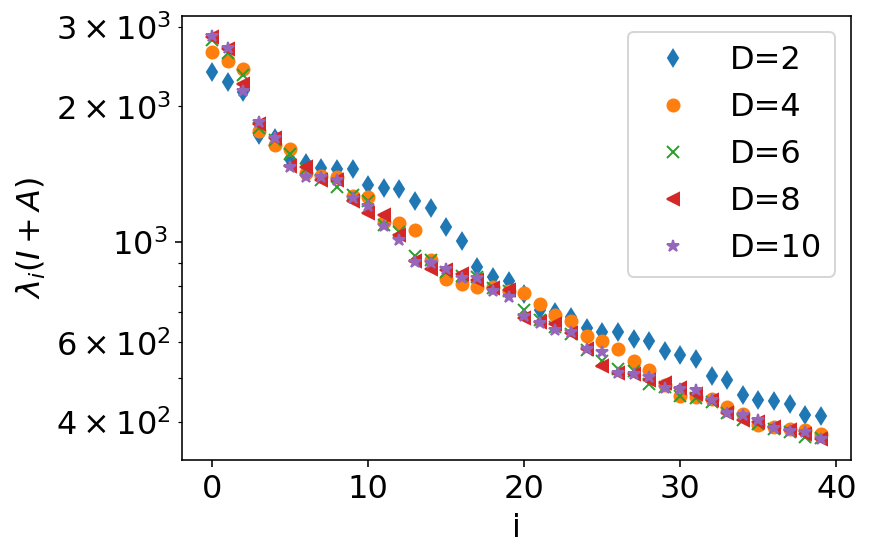}
        \caption{$M=10$; 40 eigenvalues}
    \end{subfigure}% 
    \begin{subfigure}[t]{0.5\textwidth}
        \centering
        \includegraphics[width=\linewidth]{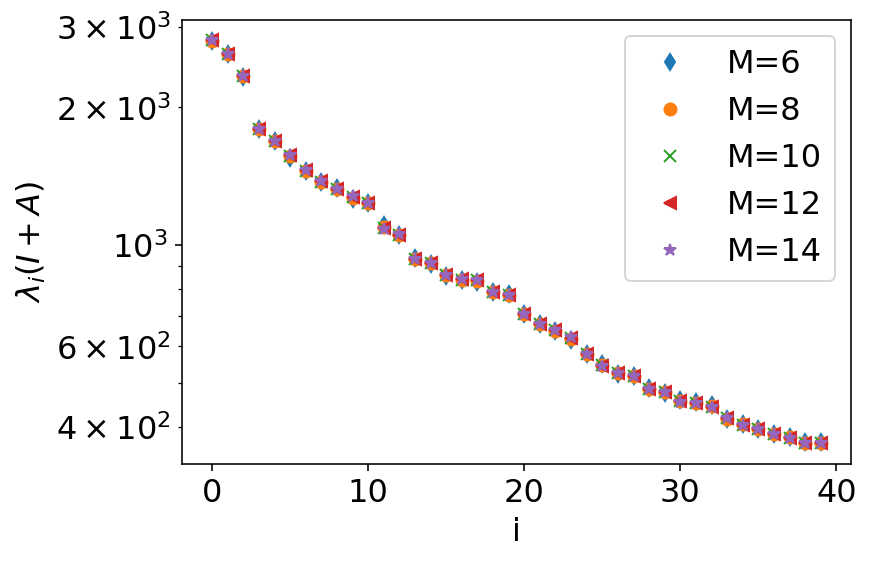}
        \caption{$D=6$; 40 eigenvalues}
    \end{subfigure}
    \caption{150 and 40 largest eigenvalues of $\bI+\bA$ with different parameter values for $\bC_\bB$. The left panels show the case for $M=10$ and changing $D$, whereas the right panels show the case with $D=6$ and changing $M$. Note that the semi-log scale is used.} \label{fig:eigs_lorenz96_Bmtx_parameters}
\end{figure}

We first explore how the parameters $D$ and $M$ for $\bC_\bB$ influence the distribution of the eigenvalues of $\bI+\bA$. These are computed via the deterministic function \texttt{scipy.sparse.linalg.eigsh}. In Figure~\ref{fig:eigs_lorenz96_Bmtx_parameters}, we show 150 largest eigenvalues of $\bI+\bA$ (all other eigenvalues are equal to one), and 40 largest eigenvalues in separate panels. It is apparent that changing the value of $M$ has a negligible effect on the distribution of the largest eigenvalues. This is not the case when $D$ is varied. In fact, larger values of $D$ result in increased gaps between eigenvalues when their decay rate changes (around $\lambda_{50}(\bI+\bA)$ to $\lambda_{60}(\bI+\bA)$). Even though we aim to compute only a relatively small number of largest eigenvalues, large trailing eigenvalues may affect the approximation quality adversely. We also note that larger values of $D$ slightly increase the largest eigenvalue and give a faster decay of the largest eigenvalues; this is a favourable setting for both randomized methods and LMPs. Further experiments use parameters $D=6$ and $M=10$.

\subsection{Results: randomized LMPs for the control member}

We now evaluate the performance of randomized spectral LMPs when solving \eqref{eq:system_with_CVT} for the control member. The scaling parameter in LMPs is set to $\theta=(\lambda_k(\bI+\bA)+1)/2$, an approximate average between the largest and smallest eigenvalues of the preconditioned system. Here and in the following sections $\lambda_i(\bI+\bA)$ denote the approximations of the eigenvalues of $\bI+\bA$ computed by Algorithm~\ref{alg:nystrom}. In Algorithm~\ref{alg:nystrom} we set $\ell = 20$ and use all the computed eigenpairs to construct the LMPs, that is, there is no oversampling and $k=20$.
The sketching matrices are as listed in section~\ref{sec:numerical_experiments}. We consider 20 random initialisations for each choice of $\bOmega$, and use the resulting randomized LMPs in PCG. The median over the random runs is reported; the mean results are qualitatively similar.

\begin{figure}
    \centering
    \includegraphics[width=0.5\linewidth]{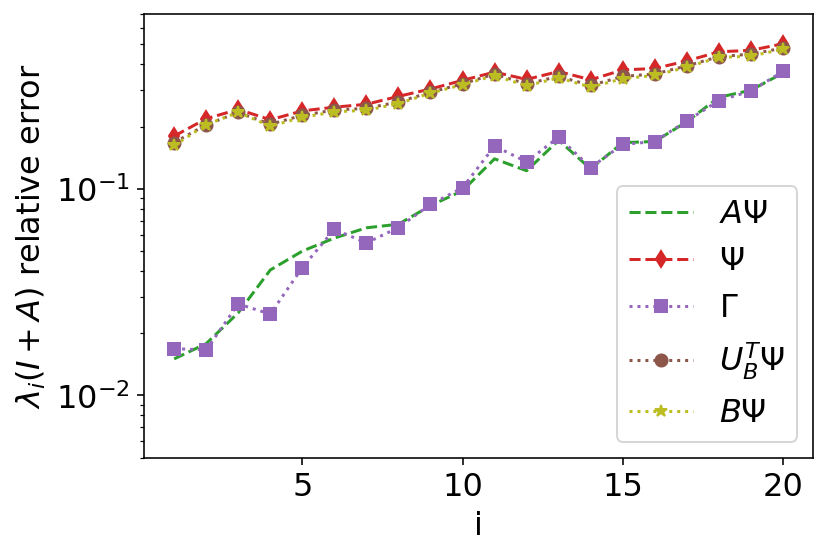}
    \caption{Median value of the relative eigenvalue approximation error with different choices of the sketching matrix $\bOmega$. }
    \label{fig:l96_eigs_approx_error_median}
\end{figure}

We first explore how the choice of the sketching matrix affects the quality of the approximation of eigenvalues of $\bI+\bA$ computed by Algorithm~\ref{alg:nystrom}. We compare $\lambda_i(\bI+\bA)$ against the eigenvalues $\lambda_i^*(\bI+\bA)$ given by the deterministic function \texttt{scipy.sparse.linalg.eigsh}. We show the median (over the 20 random initialisations) of the relative error $\vert \lambda_i^*(\bI+\bA) -  \lambda_i(\bI+\bA) \vert/\lambda_i^*(\bI+\bA)$ in Figure~\ref{fig:l96_eigs_approx_error_median}. The approximation error is smaller for the largest eigenvalues and increases as the eigenvalues decrease with all choices of $\bOmega$. Note that using $\bOmega = \bGamma$ gives approximation of a similar quality as $\bOmega = \bA\bPsi$, which is in line with our analysis in section~\ref{sec:sketching_using_eda_structure}.
These two choices outperform all other $\bOmega$, especially for the largest eigenvalues. 

\begin{figure}[h]
    \centering
    \begin{subfigure}[t]{0.49\textwidth}
        \centering
        \includegraphics[width=\linewidth]{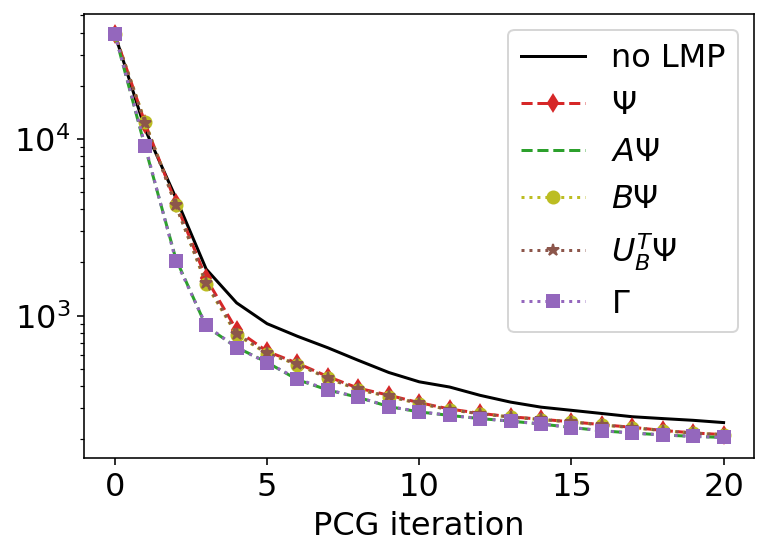}
    \end{subfigure}% 
    \begin{subfigure}[t]{0.49\textwidth}
        \centering
        \includegraphics[width=\linewidth]{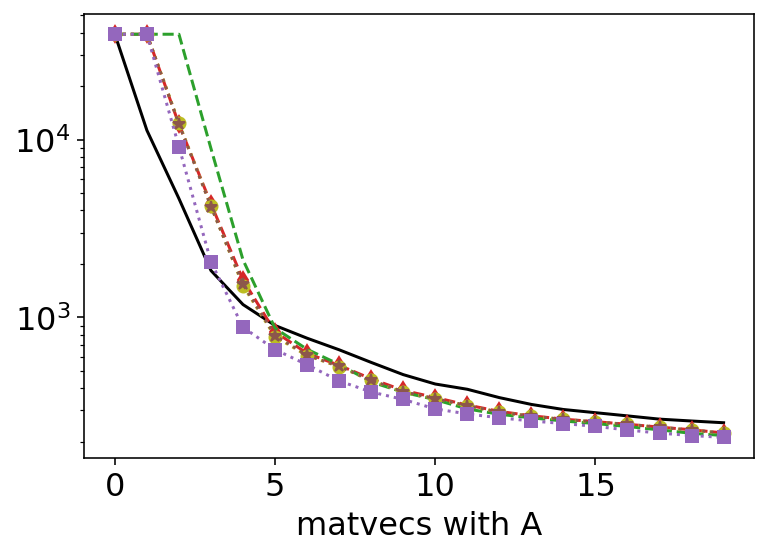}
    \end{subfigure}
\caption{The median value of quadratic cost function \eqref{eq:quadratic_cost_function} at every PCG iteration (left panel) and against matrix-vector products with $\bA$ (right panel); the median is taken over 20 random initializations of the sketching matrix $\bOmega$. The LMPs are constructed with scaling parameter $\theta=(\lambda_k(\bI+\bA)+1)/2$, and $k=20$ approximate eigenpairs returned by Algorithm~\ref{alg:nystrom}.} \label{fig:qcf_l96_k20_thmed}
\end{figure}

We now test if the differences in the eigenvalue approximation quality translate into differences in the performance of the LMPs. 
In Figure~\ref{fig:qcf_l96_k20_thmed}, we show the median value of the quadratic cost function in \eqref{eq:quadratic_cost_function} against both the number of PCG iterations and the number of matrix-vector products with $\bA$, in order to account for the cost of generating the LMPs. We assume that $\bY$ in Algorithm~\ref{alg:nystrom} is computed in parallel for each column, and thus assume that the cost of it is equivalent to computing one matrix-vector product with $\bA$. Hence step~\ref{algstep:sketch} of Algorithm~\ref{alg:nystrom} requires computing one matrix-vector product with $\bA$ regardless of the choice of $\bOmega$. Using $\bOmega = \bA \bPsi$ needs computing with $\bA$ again, whereas all other choices of $\bOmega$ do not call for an extra matrix-vector product with $\bA$; we assume that the cost of all other computations is negligible.
The randomized LMPs accelerate the reduction of \eqref{eq:quadratic_cost_function}, especially in the first iterations of PCG.
Using $\bOmega=\bGamma$ gives LMP that performs very similarly to when one iteration of the power method is run ($\bOmega=\bA\bPsi$), but at a smaller cost. Note that using sketching matrices with information coming from $\bB$ and $\bUb^\text{T}$ performs similarly to using standard Gaussian $\bOmega$ in this low-rank approximation setting. These results are consistent with the eigenvalue approximation quality results in Figure~\ref{fig:l96_eigs_approx_error_median}. We focus on $\bOmega=\bGamma$ in our further experiments.

\subsubsection{Sensitivity analysis}
We now focus on the sensitivity of the LMP when using $\bOmega=\bGamma$. Specifically, we consider how different choices of the scaling parameter $\theta$ and the use of oversampling affects the performance of LMPs. In the following, the scaling parameter in LMPs is set to one of the following values: $\theta=(\lambda_k(\bI+\bA)+1)/2$ as suggested by \cite{diouane2024efficient}, $\theta=\lambda_k(\bI+\bA)$ as suggested by \cite{diouane2024efficient} and \cite{martinsson2020randomized}, and $\theta=1$ due to already existing multiple eigenvalues equal to 1. We perform a fixed-budget experiment, where $\ell = 20$ in Algorithm~\ref{alg:nystrom}, and we use either $k=20$ (no oversampling) or $k=15$ (oversampling by 5) approximate eigenpairs. 

\begin{figure}[h]
    \centering
    \captionsetup[subfigure]{justification=centering}
    \begin{subfigure}[t]{0.5\textwidth}
        \centering
        \includegraphics[width=\linewidth]{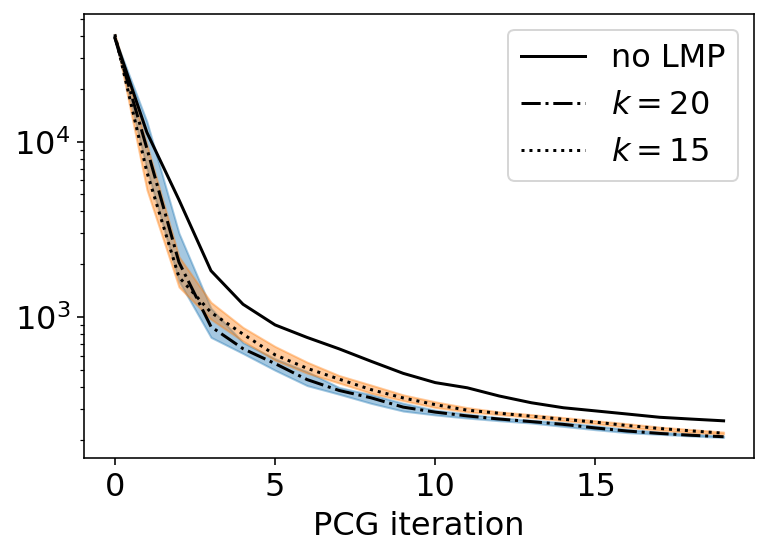}
        \caption{$\theta = (\lambda_k(\bI+\bA)+1)/2$}
    \end{subfigure}%
    \begin{subfigure}[t]{0.5\textwidth}
        \centering
        \includegraphics[width=\linewidth]{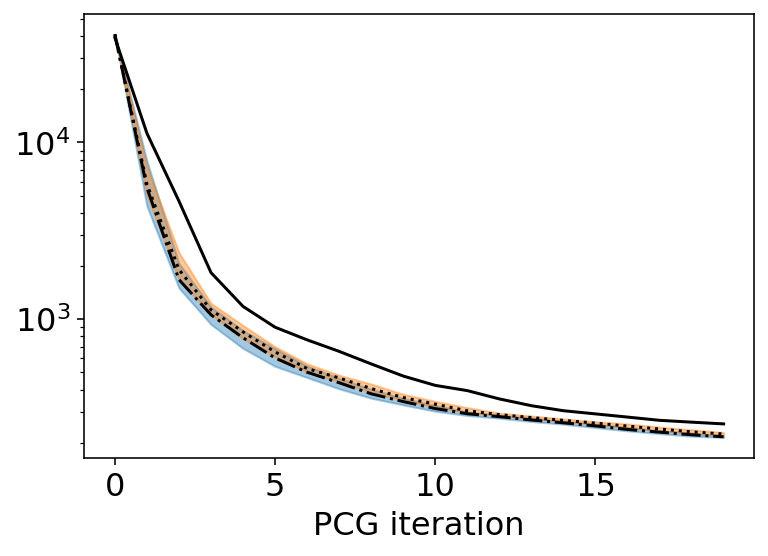}
        \caption{$\theta = \lambda_k(\bI+\bA)$}
    \end{subfigure}
     \begin{subfigure}[t]{0.5\textwidth}
        \centering
        \includegraphics[width=\linewidth]{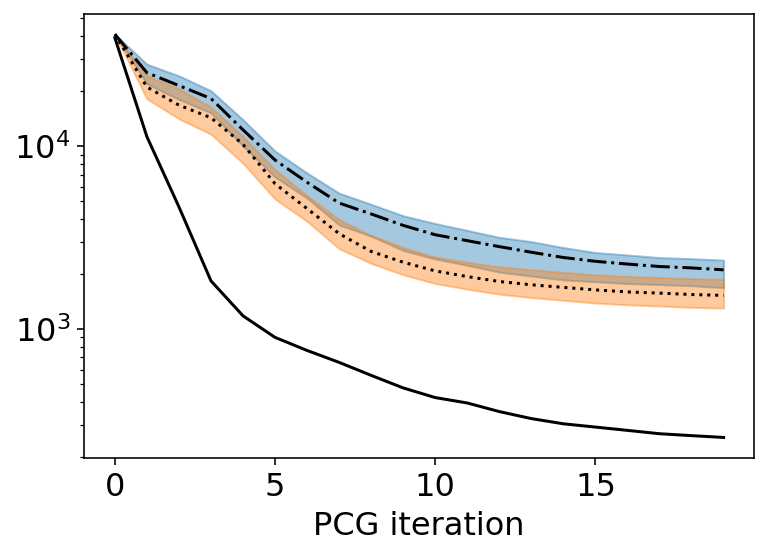}
        \caption{$\theta = 1$}
    \end{subfigure}%
    \caption{Values of quadratic cost function \eqref{eq:quadratic_cost_function} at every PCG iteration when using no LMP (solid line), and when the LMP is constructed using sketching matrix $\bOmega=\bGamma$, a different scaling parameter $\theta$ and a number of approximate eigenpairs $k$. For the cases that use LMP we show the median (dash-dotted and dotted lines), largest and smallest (coloured areas: blue for $k=\ell=20$, orange for $k=15$) values over 20 random initializations of $\bOmega$.} \label{fig:qcf_l96_k15ov5_k20_ov0_mean_max_min_thetalk_thetamed}
\end{figure}

The importance of the choice of the scaling parameter is illustrated by the LMP with $\theta=1$ performing worse than when no LMP is used in the first 20 iterations of PCG (Figure~\ref{fig:qcf_l96_k15ov5_k20_ov0_mean_max_min_thetalk_thetamed}). 
When scaling with $\theta = (\lambda_k(\bI+\bA)+1)/2$, the performance of LMP is more sensitive to the quality of the computed eigenpairs compared to $\theta=\lambda_k(\bI+\bA)$. Using $k=20$ can give worse results than using $k=15$ or no LMP in the first 2-4 PCG iterations, but afterwards using $k=20$ outperforms $k=15$. The poor initial performance when $k=20$ may be explained by the poor quality of the smallest approximate eigenvalues and associated eigenvectors.
When comparing results with $\theta = \lambda_k(\bI+\bA)$ and $\theta = (\lambda_k(\bI+\bA)+1)/2$, we notice that $\theta = \lambda_k(\bI+\bA)$ can outperform $\theta = (\lambda_k(\bI+\bA)+1)/2$ in the first few PCG iterations, but it then becomes worse in the subsequent PCG iterations. The scaling parameter thus should be tuned based on how many approximate eigenpairs can be used to construct the LMP and how many PCG iterations can be run.

\subsection{Results: randomized LMPs for the ensemble members}

\begin{figure}[h]
    \centering
    \captionsetup[subfigure]{justification=centering}
    \begin{subfigure}[t]{0.5\textwidth}
        \centering
        \includegraphics[width=\linewidth]{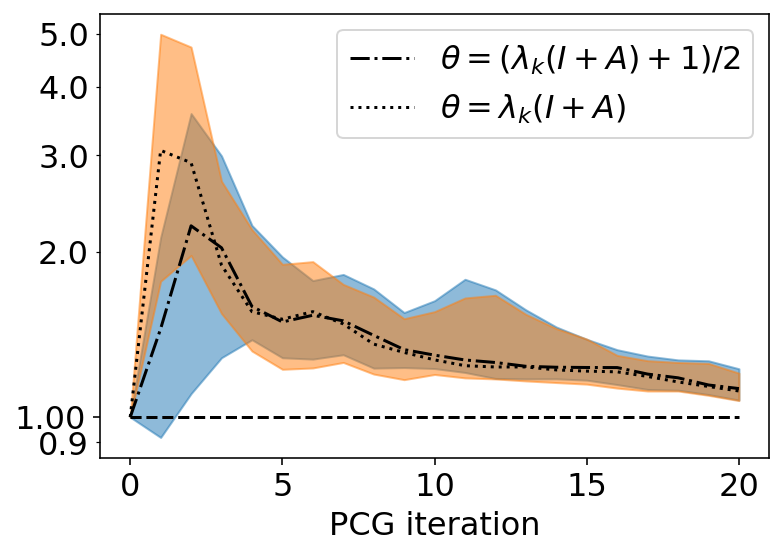}
        \caption{$k=\ell=20$}
    \end{subfigure}%
    \begin{subfigure}[t]{0.5\textwidth}
        \centering
        \includegraphics[width=\linewidth]{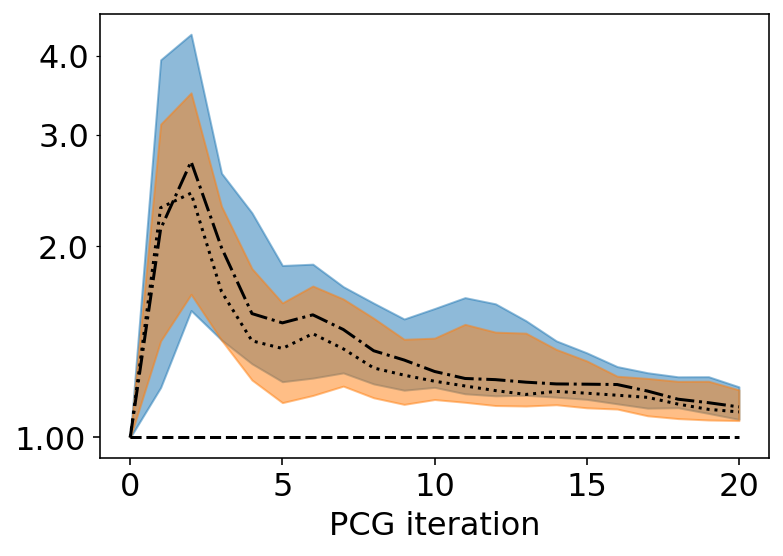}
        \caption{$k=15$}
    \end{subfigure}
    \caption{The median (dash-dotted and dotted lines), maximum and minimum values (coloured areas) over the ensemble members of the ratio between the quadratic cost functions \eqref{eq:quadratic_cost_function} at every PCG iteration when using no LMP and when using LMP constructed with $\bOmega = \bGamma$, a different scaling parameter $\theta$ and a number of approximate eigenpairs $k$. 
    The coloured areas legend: blue for $\theta = (\lambda_k(\bI+\bA)+1)/2$, orange for $\theta = \lambda_k(\bI+\bA)$. Values higher than 1 indicate larger decrease of \eqref{eq:quadratic_cost_function} when using the randomized spectral LMP compared to using no LMP. } \label{fig:qcf_l96_k15ov5_k20_ov0_mean_max_min_thetalk_thetamed_ensemble}
\end{figure}

We now test the performance of the randomized spectral LMPs when solving the ensemble problems \eqref{eq:system_with_CVT}. Note that the LMP is generated by sketching the low-rank part of the Hessian of the control problem, and is used as is for the ensemble problems. The LMPs are constructed as in the previous section using the sketching matrix of right-hand sides, with $k=20$ and $k=15$ approximate eigenpairs, and the scaling parameters $\theta = (\lambda_k(\bI+\bA)+1)/2$ and $\theta = \lambda_k(\bI+\bA)$. For one instance of the ensemble with 20 members, we compare the quadratic cost function \eqref{eq:quadratic_cost_function} values when using no LMP and when using the randomized spectral LMP (Figure~\ref{fig:qcf_l96_k15ov5_k20_ov0_mean_max_min_thetalk_thetamed_ensemble}). 
This is done by computing the ratio $\alpha_{j,i}=J\left(\dx_{j,i}^{\text{(no LMP)}}\right)/J\left(\dx_{j,i}^{\text{(LMP)}}\right)$, where $J\left(\dx_{j,i}^{\text{(no LMP)}}\right)$ is the value of \eqref{eq:quadratic_cost_function} at $i$th iteration of PCG when solving \eqref{eq:system_with_CVT} coming from the $j$th ensemble member and using no LMP; $J\left(\dx_{j,i}^{\text{(LMP)}}\right)$ corresponds to the solve that uses the LMP. At every PCG iteration $i$ we show $\max_j \alpha_{j,i}$, $\min_j \alpha_{j,i}$ and $\text{median}_j \alpha_{j,i}$.
We see that, in general, the LMP accelerates the decrease of the quadratic cost function~\eqref{eq:quadratic_cost_function} in the first PCG iterations, even though the spectral LMP has been built by sketching the Hessian matrix of the control member only. Using an LMP with $k=20$ and $\theta = (\lambda_k(\bI+\bA)+1)/2$ can have negative impact on the reduction of the quadratic cost function~\eqref{eq:quadratic_cost_function} in the first iterations of PCG. It is thus safer to use oversampling, that is, discard a few smallest approximate eigenpairs. 
The experiments have been repeated with 20 different instances of the ensemble and similar results were obtained.

\section{Conclusions}\label{sec:conclusions}
In this paper, we have considered a novel approach for accelerating the minimisation of variational data assimilation problems in an EDA setting. This approach computes a randomized approximation of the Hessian matrix and then uses it to construct a spectral LMP for the PCG solver. 

We have proposed employing a new sketching matrix within the randomized method. This sketching matrix is obtained inexpensively using the right-hand side vectors of the linear systems of equations in the EDA and contains information on the Hessian matrix. Our theoretical analysis and numerical experiments with the Lorenz-96 model have shown that the new sketching matrix can improve the accuracy of the approximation of the Hessian. This approximation is of a similar quality as the one returned by a computationally more expensive method, where the extra work is equivalent in terms of wall-clock time to approximately one PCG iteration. 
The spectral LMP constructed with the proposed randomized approximation is effective in accelerating the reduction of the quadratic cost function in the first iterations of the PCG method for both control and ensemble minimisation problems. 

Our numerical experiments have also emphasised the importance of scaling the spectral LMP and using oversampling within the randomized method. 
The performance of the LMP in the first PCG iterations is sensitive to the choices of both scaling and oversampling parameters, which should be tuned for the specific problem. In general, when solving the perturbed problems, it is safer to use oversampling. 

In this work, we have considered an EDA system which has a control member. The proposed approach can be applied to EDA systems with no reference member by replacing it with a randomly selected EDA member. However, the assumptions made in the theoretical analysis of the new sketching matrix may not hold in this case. 

The effectiveness of the proposed approach is likely to differ depending on the specific data assimilation system, the decay of the eigenvalues of the Hessian matrices and the similarity of the Hessian matrices of the EDA members. Independent ongoing experiments with a quasi-geostrophic model show impact similar to what is reported in this paper. The following stage is to validate these encouraging results in a realistic EDA setting, where the dimensions of the linear systems are significantly larger than the ones considered in this paper. 
Our ongoing work is focused on using the proposed framework for accelerating the operational EDA associated with the global data assimilation system of Météo-France.

\section*{Acknowledgements}
The Les Enveloppes Fluides et l’Environnement (LEFE) programme of the Institut National des Sciences de l’Univers—Centre National de la Recherche Scientifique (INSU-CNRS) funded part of this research.

\section*{Conflict of interest}
The authors declare no conflict of interest.

\section*{Data Availability Statement}
The code used for the simulations is available from the corresponding author upon reasonable request.

\appendix 

\section{Perturbation analysis for $\bGamma$}\label{app:perturbation_gamma}
In this section, we derive \eqref{eq:cov_gamma_Ej_term}, that is, the covariance of the RHS differences vector $\bgamma_j = \bUbj^\text{T} \bG^\text{T}_j \bR_j^{-1} \bd_j - \bU_{\bB}^\text{T} \bG^\text{T} \bR^{-1} \bd$. We relax the assumption that $\bG_j\approx\bG$ and write instead $\bG_j = \bG +\bE_j$. We still assume that $\bUbj\approx\bUb$ and $\bR_j\approx\bR$.
We use %$\bG_j = \bG +\bE_j$, 
$\bxb_j - \bxb = \beps^\text{b}_j$ and the second order approximation 
\begin{equation}
\mathcal G(\bxb_j)
\approx
\mathcal G(\bxb) + \bG \,\beps^\text{b}_j + \mathcal{S}(\beps^\text{b}_j)
\end{equation}
where
\begin{equation}
\mathcal{S}(\beps^\text{b}_j) = 
\frac12
\begin{bmatrix}
(\beps^\text{b}_j)^\text{T} \mathbf W_1 \beps^\text{b}_j \\
\vdots \\
(\beps^\text{b}_j)^\text{T} \mathbf W_m \beps^\text{b}_j
\end{bmatrix},
\end{equation}
with $\bG$ being the Jacobian of ${\mathcal G}$ evaluated at $\bxb$, and
\[
\mathbf W_i
=
\left.\nabla^2 \mathcal G_i(\bxb)\right|_{\mathbf x=\bxb}
\]
being the Hessian matrix of the $i$-th component $\mathcal{G}_i$ evaluated at $\bxb$.
Then we can write
\begin{align}
    \bgamma_j  \approx & \, \bUb^\text{T} (\bG +\bE_j)^\text{T} \bR^{-1} \bd_j - \bU_{\bB}^\text{T} \bG^\text{T} \bR^{-1} \bd\\
    \approx & \, \bUb^\text{T} (\bG +\bE_j)^\text{T} \bR^{-1} \left(\by + \bUR \boeta_j^o - (\mathcal G(\bxb) + \bG \,\beps^\text{b}_j + \mathcal{S}(\beps^\text{b}_j))\right) - \bU_{\bB}^\text{T} \bG^\text{T} \bR^{-1} \left(\by - \mathcal{G}(\bxb) \right) \\
   = & \,\bUb^\text{T} \bG^\text{T} \bR^{-1} \left( \bUR \boeta^o_j - \bG \beps^\text{b}_j - \mathcal{S}(\beps^\text{b}_j)\right)
                    + \bUb^\text{T} \bE_j^\text{T} \bR^{-1} \left(\by + \bUR \boeta^o_j - \mathcal{G}(\bxb) - \bG \beps^\text{b}_j - \mathcal{S}(\beps^\text{b}_j)\right) \\
              \approx  & \, \bUb^\text{T} \bG^\text{T} \bR^{-1} \left( \bUR \boeta^o_j - \bG \beps^\text{b}_j - \mathcal{S}(\beps^\text{b}_j)\right)
                    + \bUb^\text{T} \bE_j^\text{T} \bR^{-1} \left(\by + \bUR \boeta^o_j - \mathcal{G}(\bxb) - \bG \beps^\text{b}_j \right), 
\end{align}
where we omit the third-order term $\bUb^\text{T} \bE_j^\text{T} \bR^{-1} \mathcal{S}(\beps^\text{b}_j)$ in the last approximation. 
Recalling the assumptions $\mathbb{E}[\boeta^o_j]=0$ and $\mathbb{E}[\beps^\text{b}_j]=0$, and the independence of $\boeta^o_j$ and $\bE_j$, the expected value is 
\begin{equation}
   \mathbb{E}\! \left[ \bgamma_j\right]  \approx  -  \bUb^\text{T} \bG^\text{T} \bR^{-1}\mathbb{E}\! \left[ \mathcal{S}(\beps^\text{b}_j)\right] + \bUb^\text{T} \mathbb{E}\! \left[ \bE_j^\text{T} \bR^{-1} \left(\bd - \bG \beps^\text{b}_j \right) \right].
\end{equation}
%%%%%%%%%%%%%%%%%%%%%%%%%%%%%%%%%%%%%%%%%%%%%%%%%%%%%%%%%%%%%%%%%%%%%%%%%%%%%%%%%%%%%%%%%%%%%%%%%%%%%%%%%%%%%%%%%%%%%%%%%%%%%%%%%%%%%%%%%%%%%%%%%%%%%%%%%%%%%%%%%%%%%%%%%%%%%%

Using $\operatorname{cov}(\beps_j^\text{b},\beps_j^\text{b}) \approx \bB$ and $\boeta^o_j$ being independent of both $\bxb_j$ and $\bE_j$, we then compute the covariance as
%\begin{equation*}
$    \operatorname{cov}(\bgamma_j,\bgamma_j) =\mathbb{E}\! \left[\left( \bgamma_j -  \mathbb{E}\![\bgamma_j]\right) \left( \bgamma_j -  \mathbb{E}\![\bgamma_j]\right)^\text{T} \right]$,
%\end{equation*}
term-by-term, and obtain
\begin{align}
    \operatorname{cov}(\bgamma_j,\bgamma_j) & \approx  \bA + \bA^2 + \bDelta_1 + \bDelta_1^\text{T}  + \bDelta_2,
\end{align}
where
\begin{align}
    \bDelta_1  =  & \ \bUb^\text{T} \bG^\text{T} \bR^{-1}\mathbb{E}\! \left[ \bE_j  \right]\bUb  
    - \bUb^\text{T} \bG^\text{T} \bR^{-1}\bG\mathbb{E}\! \left[  \beps^\text{b}_j \bd^\text{T} \bR^{-1} \bE_j \right] \bUb 
    + \bUb^\text{T} \bG^\text{T} \bR^{-1} \bG\mathbb{E}\! \left[ \beps^\text{b}_j (\beps^\text{b}_j)^\text{T} \bG^\text{T} \bR^{-1} \bE_j \right] \bUb \nonumber \\
    & - \bUb^\text{T}\mathbb{E}\! \left[  \bE_j^\text{T} \bR^{-1} \bG \beps^\text{b}_j \bd^\text{T} \bR^{-1} \bE_j \right] \bUb 
    + \bUb^\text{T}\mathbb{E}\! \left[   \bE_j^\text{T}  \right] \bR^{-1} \bd
   \mathbb{E}\! \left[ (\beps^\text{b}_j)^\text{T} \bG^\text{T} \bR^{-1} \bE_j \right] \bUb \nonumber \\
    & + \bUb^\text{T}\bG^\text{T}\bR^{-1}\mathbb{E}\!\left[\mathcal{S}(\beps^\text{b}_j)\right]\bd^\text{T}\bR^{-1}\mathbb{E}\!\left[\bE_j\right] \bUb - \bUb^\text{T}\bG^\text{T}\bR^{-1}\mathbb{E}\!\left[\mathcal{S}(\beps^\text{b}_j)\bd^\text{T}\bR^{-1}\bE_j\right]\bUb  \nonumber \\
    & +\bUb^\text{T}\bG^\text{T}\bR^{-1}\mathbb{E}\!\left[\mathcal{S}(\beps^\text{b}_j) (\beps^\text{b}_j)^\text{T}\right] \bG^\text{T} \bR^{-1}\bG \bUb,\\
%\end{align*}
%\begin{align*}
    \bDelta_2 =  & \ \bUb^\text{T}\mathbb{E}\! \left[  \bE_j^\text{T} \bR^{-1} \bE_j \right] \bUb 
    + \bUb^\text{T}\mathbb{E}\! \left[  \bE_j^\text{T} \bR^{-1} \bd \bd^\text{T} \bR^{-1} \bE_j \right] \bUb - \bUb^\text{T}\mathbb{E}\!\left[  \bE_j^\text{T}  \right] \bR^{-1} \bd \bd^\text{T} \bR^{-1}\mathbb{E}\! \left[ \bE_j\right] \bUb 
\end{align}
and where only terms up to the third-order in $\beps^\text{b}_j$ and $\bE_j$ are kept.

\section{Nystr\"{o}m algorithms}\label{app:Nystrom_algorithms}
A few versions of algorithms that return the Nystr\"{o}m approximation in \eqref{eq:nystrom_approx} can be found in the literature. We highlight some of the main differences in the following.
\begin{itemize}
    \item The sketching matrix $\bOmega$ can be chosen such that it samples $\ell$ columns of $\bA$. This approach is common when dealing with kernel matrices, for example, \citep{williams2000using,drineas2005nystrom}. Note that we consider a general random projection approach in this paper, where $\bOmega$ involves a random Gaussian matrix. 
    \item The matrix $\bA$ can be accessed just once giving a so-called single-pass algorithm, or more times if using the power method. Note that in the latter case or if $\bOmega$ has a very large number of columns, then orthogonalizing the columns of $\bA \bOmega$ as described in Algorithm~\ref{alg:nystrom} may be necessary for numerical stability, e.g. Chapter 6 of \cite{stewart2001matrix}, \cite{tropp2017fixed}.
    \item The small matrix $\bOmega^\text{T} \bA \bOmega$ may be inverted via different factorizations. Cholesky decomposition is often used (as in Algorithm~\ref{alg:nystrom}), e.g. \citep{halko2011finding}, but truncated eigendecomposition is also considered \citep{persson2023randomized}.
    \item \cite{li2017algorithm} proposed working with shifted $\bA +\nu \bI $, $\nu >0$, instead of $\bA$ to ensure that the Cholesky decomposition does not break down. Some suggested shifts are $\nu = \epsilon \Vert \bA \bOmega \Vert_F$ and $\nu = \epsilon \operatorname{trace}(\bA)$, where $\epsilon$ is the machine epsilon \citep{tropp2023randomized,tropp2017fixed}. In our experiments, where the sample size $\ell$ is very small, the eigenvalues of $\bOmega^\text{T} \bA \bOmega$ are always larger than one, thus the shift is not needed.
\end{itemize}

\printendnotes

\bibliography{sample}
\graphicalabstract{plots_lorenz96/ensemble/qcf_ratio_ensemble_k15_rnd1.png}{Randomized matrix approximation methods are used for preconditioning iterative solvers in ensembles of data assimilations (EDA), exploiting its specific random structure.
Improved convergence of iterative solvers for all EDA members is demonstrated numerically with the Lorenz-96 model when using this randomisation-based limited memory preconditioner.}

\end{document}